\documentclass{elsart}

\usepackage{amsmath}
\usepackage{amssymb}
\usepackage{amscd}
\usepackage{sty1}

\journal{Journal of Functional Analysis}

\makeatletter
\@addtoreset{equation}{section}
\makeatother

\def\bloneg{\mathrm{L}^1(G)}

\def\bltwog{\mathrm{L}^2(G)}
\def\wcbop#1{\mathcal{CB}^\sigma(#1)}
\def\cbnorm#1{\norm{#1}_{cb}}
\def\cring{\Omega(G)}
\def\cringftwo{\Omega(\ftwo)}
\def\cringh{\Omega(H)}
\def\cringhcg{\Omega(H\cross G)}
\def\cstarg{\mathrm{C}^*(G)}

\def\disunion#1#2{\;\cdot\!\!\!\!\!\!\bigcup_{#1}^{#2}}
\def\falftwo{\mathrm{A}(\ftwo)}
\def\fal#1{\mathrm{A}(#1)}
\def\falcg{\mathrm{A}_c(G)}
\def\falch{\mathrm{A}_c(H)}
\def\falg{\mathrm{A}(G)}
\def\falh{\mathrm{A}(H)}
\def\falgg{\mathrm{A}(G\cross G)}
\def\fsal#1{\mathrm{B}(#1)}
\def\fsalg{\mathrm{B}(G)}
\def\fsalhcg{\mathrm{B}(H\cross G)}
\def\fsalhcgd{\mathrm{B}(H_d\cross G_d)}
\def\fsalh{\mathrm{B}(H)}
\def\fsaly{\mathrm{B}(Y)}
\def\ftwo{\mathrm{F}_2}
\def\id{\mathrm{id}}
\def\smatn{\mathrm{M}_n}
\def\matn#1{\mathrm{M}_n(#1)}
\def\meash{\mathrm{M}(H)}
\def\ocring{\Omega_o(G)}
\def\ocringh{\Omega_o(H)}
\def\ptens{\widehat{\otimes}}
\def\tens{\otimes}
\def\trace{\mathrm{tr}}
\def\vnftwo{\mathrm{VN}(\ftwo)}
\def\vng{\mathrm{VN}(G)}
\def\vngg{\mathrm{VN}(G\cross G)}
\def\vnh{\mathrm{VN}(H)}
\def\vnpi{\mathrm{VN}_\pi}
\def\vntens{\overline{\otimes}}
\def\wapg{\mathrm{W}\mathrm{A}\mathrm{P}(G)}
\def\wstarg{\mathrm{W}^*(G)}
\def\wstarh{\mathrm{W}^*(H)}
\def\wstary{\mathrm{W}^*(Y)}
\def\wstarhcg{\mathrm{W}^*(H\cross G)}

\newtheorem{coset}{Proposition}[section]
\newtheorem{graphs}[coset]{Lemma}
\newtheorem{context}[coset]{Lemma}

\newtheorem{contractiveidempotent}{Theorem}[section]
\newtheorem{fcfail}[contractiveidempotent]{Proposition}

\newtheorem{affinecbhomo}{Proposition}[section]
\newtheorem{affinecbhomo1}[affinecbhomo]{Corollary}
\newtheorem{notnormext}[affinecbhomo]{Remark}
\newtheorem{bad}[affinecbhomo]{Lemma}
\newtheorem{cpmap}[affinecbhomo]{Lemma}
\newtheorem{inversemap}[affinecbhomo]{Remark}
\newtheorem{maintheorem}[affinecbhomo]{Theorem}
\newtheorem{failure}[affinecbhomo]{Proposition}
\newtheorem{maintheorem1}[affinecbhomo]{Corollary}
\newtheorem{hconnected}[affinecbhomo]{Corollary}
\newtheorem{mapfalg}[affinecbhomo]{Corollary}
\newtheorem{walter}[affinecbhomo]{Corollary}
\newtheorem{image}[affinecbhomo]{Theorem}

\begin{document}
\begin{frontmatter}

\title{Completely Bounded Homomorphisms of the Fourier Algebras}

\author[lake]{Monica Ilie\thanksref{am}} and
\ead{ilie@math.tamu.edu}
\author[wat]{Nico Spronk\thanksref{am}\corauthref{cor}}
\corauth[cor]{Corresponding author.}
\ead{nspronk@math.uwaterloo.ca}

\thanks[am]{Both authors were supported by NSERC PDFs
and were Visiting Assitant Professors
at Texas A\&M University when this work was completed.}

\address[lake]{Department of Mathematical Sciences, Lakehead University, 
Thunder Bay ON, P7B 5E1, Canada}
\address[wat]{Department of Pure Mathematics, University of Waterloo, Waterloo 
ON, \\ N2L 3G1, Canada}






\begin{abstract}
For locally compact groups $G$ and $H$ let $\falg$ denote the Fourier algebra
of $G$ and $\fsalh$ the Fourier-Stieltjes algebra of $H$.  Any continuous
piecewise affine map $\alp:Y\subset H\to G$ (where $Y$ is an element of
the open coset ring) induces a completely bounded homomorphism
$\Phi_\alp:\falg\to\fsalh$ by setting $\Phi_\alp u=u\comp\alp$ on $Y$
and $\Phi_\alp u=0$ off of $Y$.  We show that if $G$ is amenable
then any completely bounded homomorphism $\Phi:\falg\to\fsalh$ is of this
form; and this theorem fails if $G$ contains a discrete nonabelian free group.
Our result generalises results of P.J.\ Cohen (1960, {\it Amer.\ J.\ Math.},
82:213-226), B.\ Host (1986, {\it Bull.\ Soc.\ Math.\ France}, 114)
and of the first author (2004, {\it J.\ Funct.\ Anal.}, 213).  
We also obtain a description of all the idempotents
in the Fourier-Stieltjes algebras which are contractive or positive definite.
\end{abstract}


\end{frontmatter}

\pagebreak

\section{Introduction}

For any pair of locally compact abelian groups $G$ and $H$, Cohen \cite{cohen}
characterised all of the bounded homomorphisms from the group algebra
$\bloneg$ to the measure algebra $\meash$ in terms of piecewise affine maps
between their dual groups $\what{H}$ and $\what{G}$.  In doing so he made
use of an equally profound discovery of his \cite{cohen1} characterising
idempotent measures on abelian groups, a result which won him the B\^{o}chner
Memorial Prize in 1964.  These results generalised results of many
authors \cite{rudin1,rudin2,beurlingh,kahane,kahane1,leibenson}.  
Cohen's work is also exposed nicely in \cite{rudin}.   As this
characterisation is in terms of dual groups, it
is more naturally formulated in terms of the algebras of Fourier and
Fourier-Stieltjes transforms $\fal{\what{G}}$ and $\fsal{\what{H}}$.
There is a formulation of the Fourier and Fourier-Stieltjes algebras
$\falg$ and $\fsalh$, due to Eymard \cite{eymard}, which can be done for
any locally compact groups, and which generalises the group algebra and the 
measure algebra for dual groups.  Moreover, these are commutative semi-simple 
Banach algebras, even for non commutative groups.  A longstanding question
in harmonic analysis has been to determine to which extent Cohen's theorem
holds in the nonabelian setting.  Various intermediate results, described
below, have been given over the years, and our main objective in this paper
is to give a definitive solution to this problem.

The first decisive step in generalising Cohen's results
is due to Host \cite{host}.  He
first discovered the general form of idempotents in the Fourier-Stieltjes
algebras, a significant result.
He then identified the role of tensor products in obtaining
the characterisation of bounded homomorphisms from $\falg$ to $\fsalh$.
Unfortunately, by a result of Losert \cite{losert}, the Banach algebra
$\falg\tens^\gamma\falg$ (projective tensor product) is isomorphic
to $\falgg$ only if $G$ has an abelian subgroup of finite index.  Hence
it was only for such groups that Host established his result on homomorphisms.

In the intervening years the theory of operator spaces and completely
bounded maps was developed by Paulsen, Blecher, Effros and Ruan amongst 
many others.  See \cite{effrosrB}.  It was recognised in 
\cite{blecher,effrosr1} that $\falg$ and $\fsalh$ can be regarded as operator 
spaces.  The first major application of this was given by Ruan in \cite{ruan}
where it was shown that the operator space structure on $\falg$ gives rise to
a more tractable cohomology theory than traditional Banach algebra cohomology
(such as in \cite{johnson}).   
It thus makes sense to speak about {\it completely bounded
homomorphisms} from $\falg$ to $\fsalh$.  Any results on such generalise
the results of Cohen and Host since it was shown by Forrest and Wood
\cite{forrestw} that any bounded linear map from $\falg$ to any operator space,
is automatically completely bounded if and only
if $G$ has an abelian subgroup of finite index.  The advantage of
having the context of operator spaces is that it gives us the operator
projective tensor product $\ptens$ \cite{blecherp,effrosr1}.
By this, $\falg\ptens\falg$ can be naturally identified with
$\falgg$ \cite{effrosr}.  Using these techniques, 
the first author \cite{ilie,ilieTh}
took the next decisive step and characterised all of the completely bounded 
homomorphisms form $\falg$ to $\fsalh$, when $G$ is an amenable discrete
group.

In this article we generalise this result to what appears to be the 
fullest extent possible.  We make note of the fact that for any
pair of locally compact groups $G$ and $H$, a continuous piecewise affine
map $\alp:Y\subset H\to G$ induces a completely bounded homomorphism
from $\falg$ to $\fsalh$
(Proposition \ref{prop:affinecbhomo}).  We then show that if $G$ is amenable
then every  completely bounded homomorphism from $\falg$ to $\fsalh$
is thus induced (Main Result: Theorem \ref{theo:maintheorem}).  
Along the way we make significant use of a positive {\it bounded approximate
diagonal} for $\falg\ptens\falg$, following the construction of
Aristov, Runde and the second author \cite{aristovrs} 
(Lemma \ref{lem:bad}) -- which simplifies a construction of Ruan 
\cite{ruan}.  As a complement to our main theorem we
show that for any group which contains a discrete nonabelian free group
the main result fails (Proposition \ref{prop:failure}), lending
strong evidence that amenability is an indispensable assumption.  
This makes use of Leinert's free sets \cite{leinert}.
We also indicate how, for most (amenable) groups, our main result can fail
for bounded homomorphisms from $\falg$ to itself which are not
completely bounded (Remark \ref{rem:inversemap}).

In order to refine
our main result, i.e.\ to characterise completely contractive and 
``completely positive''
homomorphisms, we obtain a description of contractive and of
positive definite idempotents in Fourier-Stieltjes algebras (Theorem
\ref{theo:contractiveidempotent}).  This result is well known for abelian
groups, but does not appear to be in the literature for general groups.
It is a special case of the significant theorem of
Host \cite{host}, though it is not mentioned nor covered by him.

\subsection{Preliminaries}
If $G$ is any locally compact group let $\falg$ denote its {\it Fourier
algebra} and $\fsalg$ denote its {\it Fourier-Stieltjes algebra}, as
defined in \cite{eymard}.  We recall that $\fsalg$ consists of all
matrix coefficients of continuous unitary representations, i.e.\
functions of the form $s\mapsto\inprod{\pi(s)\xi}{\eta}$ where
$\pi:G\to\fU(\fH)$ is a homomorphism, continuous when the unitary group
$\fU(\fH)$ on the Hilbert space $\fH$ is endowed with the weak operator 
topology.  We also recall that $\falg$ is the space of all matrix
coefficients of the left regular representation $\lam_G:G\to\fU(\bltwog)$,
given by left translation operators on
$\bltwog$, the Hilbert space of (equivalence classes of)
square-integrable functions.  The norms on $\falg$ and $\fsalg$ are given
by the dualities indicated below.

We note that $\falg$ has bounded dual space
$\falg^*\cong\vng$, where $\vng$ is the von Neumann algebra generated by 
$\lam_G$.  The Fourier-Stieltjes algebra
is the predual of the {\it enveloping von Neumann algebra} $\wstarg$ which
is generated by the universal representation $\varpi_G$ \cite{dixmierB}.
On the other hand, $\fsalg$ is the dual of the  {\it enveloping C*-algebra}
$\cstarg$.  We note that $\wstarg$ satisfies the universal property for group
von Neumann algebras:  if $\pi:G\to\fU(\fH)$ is a continuous 
representation and $\vnpi$ is the von Neumann algebra it generates, then
there is $*$-homomorphism $\Pi:\wstarg\to\vnpi$ such that
$\Pi(\varpi(s))=\pi(s)$ for each $s\iin G$.  We note that both
$\falg$ and $\fsalg$ are semi-simple commutative Banach algebras under
pointwise operations.  Moreover, $\falg$ is an ideal in $\fsalg$. 
Furthermore, $\falg$ has 
Gel'fand spectrum $G$, and is {\it regular} on $G$ in the sense that
for any compact subset $K$ of $G$, and any open set $U$ containing $K$, there
is an element $u\iin\falg$ such that $u|_K=1$ and $\supp{u}\subset U$.

Our standard reference for operator spaces and completely bounded maps
is \cite{effrosrB}, though we will indicate other references below.  
Any C*-algebra $\fA$ is an {\it operator space}
in the sense that for $n=1,2,\dots$ the algebra of $n\cross n$ matrices
over $\fA$, $\matn{\fA}$ admits a unique norm which makes it into a 
C*-algebra.  A linear map $T:\fA\to\fB$ between C*-algebras is called
{\it completely bounded} if it is bounded and 
its amplifications $T^{(n)}:\matn{\fA}\to
\matn{\fB}$, given by $T^{(n)}[a_{ij}]=[Ta_{ij}]$, give a bounded
family of norms $\left\{\norm{T^{(n)}}:n=1,2,\dots\right\}$.  In this case
we write $\cbnorm{T}=\sup\left\{\norm{T^{(n)}}:n=1,2,\dots\right\}$.
Moreover we say $T$ is {\it completely contractive} if $\cbnorm{T}\leq 1$;
that $T$ is a {\it complete isometry} if each $T^{(n)}$ is an isometry;
and that $T$ is {\it completely positive} if each $T^{(n)}$ is a positive
map (see \cite{paulsenB}).  Examples of completely bounded maps are
$*$-homomorphisms, which are also completely positive and contractive,
and multiplications by fixed elements in C*-algebras.

If $\fM$ and $\fN$ are von Neumann algebras with
preduals $\fM_*$ and $\fN_*$, we say a map $\Phi:\fM_*\to\fN_*$ is
{\it completely bounded (contractive)}, if its adjoint 
$\Phi^*:\fN\to\fM$ is such.  However, it is often convenient to
consider completely bounded maps on the space $\fM_*$ by noting that
it admits an operator space structure via the identifications
$\matn{\fM_*}\cong\wcbop{\fM,\smatn}$, $n=1,2,\dots$, where 
$\wcbop{\fM,\smatn}$ is the space of normal completely bounded maps form 
$\fM$ to the finite dimensional von Neumann algebra of $n\cross n$
complex matrices \cite{blecher,effrosr1}.  We note that these
spaces $\fM_*$, with the above matricial structures,
are completely isometrically isomorphic to subspaces of C*-algebras 
\cite{ruan0}, which are not generally operator subalgebras.
In particular, these
structures are used to create the {\it operator projective tensor product}
$\fM_*\ptens\fN_*$ \cite{blecherp,effrosr1}.  This tensor product
admits the very useful formula $(\fM_*\ptens\fN_*)^*\cong\fM\vntens\fN$
\cite{effrosr}, where $\fM\vntens\fN$ is the von Neumann tensor product.
We note for any locally compact group $G$ that multiplication
extends to a completely contractive linear map 
$\mu:\fsalg\ptens\fsalg\to\fsalg$.  Indeed 
$\mu^*:\wstarg\to\wstarg\vntens\wstarg$ is the $*$-homomorphism which extends
$\varpi_{G\times G}(s,t)\mapsto\varpi_G(s)\tens\varpi_G(t)$.
Hence we say that $\fsalg$ is a {\it completely contractive Banach algebra}. 
In particular,
multiplication by a fixed element $v\mapsto uv$ on $\fsalg$ is completely
bounded with $\cbnorm{v\mapsto uv}=\norm{u}_{\fsalg}$.

\subsection{Piecewise Affine Maps}

In this section we give a quick survey of piecewise affine maps on groups.
These maps are natural generalisations of group homomorphisms and
are the natural morphisms on finite collections of cosets.  We will require
general versions of several results from \cite{rudin} concerning abelian
groups.  Unfortunately modification of the original proofs is required,
and we give these below.

Let $G$ be a group.  A {\it coset} of $G$ is any subset $C$ of $G$ for which
there is a subgroup $H$ of $G$, and an element $s\iin G$ such that
$C=sH$.  We note that for $H$ and $s$ as above, we have that
$Hs=ss^{-1}Hs$, which means that we need not distinguish between left and right
cosets.  The following result is \cite[3.7.1]{rudin} in the case that $G$ is
abelian.  

\begin{coset}\label{prop:coset}
A subset $C$ of $G$ is a coset if and only if for every $r,s$ and $t\iin C$,
$rs^{-1}t\in C$ too.  Moreover, $C^{-1}C$ is a subgroup for which
$C=sC^{-1}C$ for any $s\iin C$.
\end{coset}

\proof Necessity is trivial, so we will prove only sufficiency.  We will
show that $H=C^{-1}C$ is a subgroup and $C=sH$ for any $s\iin C$.  If
$s,t\in H$, then $s=s_1^{-1}s_2$ and $t=t_1^{-1}t_2$ where $s_i,t_i\in C$
for $i=1,2$.  Then
\[
st=s_1^{-1}(s_2t_1^{-1}t_2)\in C^{-1}C\quad\aand\quad
s^{-1}=s_2^{-1}s_1\in C^{-1}C
\]
whence $H=C^{-1}C$ is a subgroup.  Now if $s\in C$ and $t\in H$ with
$t=t_1^{-1}t_2$ as above, then $st=st_1^{-1}t_2\in C$, so $sH\subset C$.
Also, $C=ss^{-1}C\subset sC^{-1}C=sH$.  Hence $C=sH$. \endpf

Now let $H$ be another group.  A map $\alp:C\subset H\to G$ is called
{\it affine} if $C$ is a coset and for $r,s,t\iin C$
\[
\alp(rs^{-1}t)=\alp(r)\alp(s)^{-1}\alp(t).
\]
It is clear from Proposition \ref{prop:coset} above, that the range $\alp(C)$
of $\alp$ is also a coset.  Hence if $s\in C$, then
\begin{equation}\label{eq:affhomo}
s^{-1}C\ni t\mapsto \alp(s)^{-1}\alp(st)\in \alp(s)^{-1}\alp(C)
\end{equation}
is a homomorphism between subgroups.

We let $\cringh$ denote the {\it coset ring} of the group $H$, which is
the smallest ring of subsets which contains every coset.  A map
$\alp:Y\subset H\to G$ is called {\it piecewise affine} if
\begin{align}
\text{(i) } & \text{there are pairwise disjoint }Y_i\in\cringh, 
\ffor i=1,\dots,n \notag \\
&\text{such that }Y=\disunion{i=1}{n}Y_i\text{ (disjoint union)}, and
\label{eq:piecwiseaff} \\
\text{(ii) } & \text{each }Y_i\text{ is contained in a coset }L_i
\text{ on which there is an} \notag \\
&\text{affine map } \alp_i:L_i\to G \text{ such that }
\alp_i|_{Y_i}=\alp|_{Y_i}.
\notag
\end{align}
If $\alp:Y\subset H\to G$ is a function we define the {\it graph}
of $\alp$ to be the set
\begin{equation}\label{eq:graph}
\Gamma_\alp=\{(s,\alp(s)):s\in Y\}.
\end{equation}
The following lemma is given for abelian groups in \cite[4.3.1]{rudin}.
Our proof is adapted from the one given there.

\begin{graphs}\label{lem:graphs}
Let $\alp:Y\subset H\to G$ be a function.  Then $\alp$ enjoys the following
properties.

{\bf (i)} If $\Gamma_\alp$ is a subgroup then so too is $Y$ and
$\alp$ is a homomorphism of subgroups.

{\bf (ii)} If $\Gamma_\alp$ is a coset then so too is $Y$ and
$\alp$ is an affine map.

{\bf (iii)} If $\Gamma_\alp\in\cringhcg$, then $\alp$ is 
a piecewise affine map. 
\end{graphs}

\proof (ii) Let $r,s,t\in Y$.  Then since $\Gamma_\alp$ is a coset,
\[
(r,\alp(r))(s,\alp(s))^{-1}(t,\alp(t))
=(rs^{-1}t,\alp(r)\alp(s)^{-1}\alp(t))\in\Gamma_\alp
\]
which implies that $Y$ is coset since $rs^{-1}t\in Y$.  Since $\Gamma_\alp$
is a graph, \linebreak
$(rs^{-1}t,\alp(rs^{-1}t))\in\Gamma_\alp$ too and
$\alp(r)\alp(s)^{-1}\alp(t)=\alp(rs^{-1}t)$.

(i) If $\Gamma_\alp$ is a subgroup, it is a coset containing the identity,
whence so too is $Y$.  It follows that $\alp(e)=e$ 
and $\alp$ is a homomorphism.

(iii)  Since $\Gamma_\alp\in\cringhcg$, there exists a finite collection
of subgroups $\Sigma$ of $H\cross G$ such that $\Gamma_\alp\in\fR(\Sigma)$,
the smallest ring of subsets generated by cosets of elements of $\Sigma$.
We may assume that $\Sigma$ is closed under intersections. We may also assume
that if one element of $\Sigma$ is a subgroup of another, then the index
of the first subgroup in the second is infinite.

It is then possible to write
\[
\Gamma_\alp=\disunion{i=1}{n} E_i\quad\text{where each}\quad
E_i=L_i\setdif\bigcup_{j=1}^{m_i}M_{ij}
\]
and each $L_i$ and $M_{ij}$ are cosets of elements of $\Sigma$ with
$M_{ij}\subset L_i$ for each $i,j$.  Note that each $E_i$ is itself a graph.

We claim that each $L_i$ is a graph.  If not, there are elements
$(s,t_1)$ and $(s,t_2)$ in $L_i$ with $t_1\not=t_2$, so $(e,t)
=(e,t_1^{-1}t_2)\in L_i^{-1}L_i$.  Now if $(s,\alp(s))$ is any
element of $E_i$, then $(s,\alp(s)t)\in L_iL_i^{-1}L_i=L_i$, so
$(s,\alp(s)t)\in M_{ij}$ for some $j$, since $E_i$ is a graph.
Hence 
\[
(s,\alp(s))=(s,\alp(s)t)(e,t)^{-1}\in M_{ij}(e,t)^{-1}.
\]
We note that $M_{ij}(e,t)^{-1}$ may not be in $\fR(\Sigma)$, but it
is a coset of a subgroup of infinite index in $L_i^{-1}L_i$.  Thus
\[
E_i\subset\bigcup_{j=1}^{m_i}M_{ij}(e,t)^{-1}
\quad\text{whence}\quad
L_i\subset\bigcup_{j=1}^{m_i}\brac{M_{ij}\cup M_{ij}(e,t)^{-1}}
\]
Hence the subgroup $L_i^{-1}L_i$ can be covered by finitely
many cosets of subgroups which are of infinite index in itself, 
which is impossible by
\cite{neumann} (also see Proposition \ref{prop:fcfail} for an
analytic proof of this).  Thus $L_i$ is a graph and we may write
\[
L_i=\{(s,\alp_i(s)):s\in K_i=p_1(L_i)\}
\]
where $p_1:H\cross G\to H$ is the standard projection.  Now if we let
$Y_i=p_1(E_i)$ and $N_{ij}=p_1(M_{ij})$ we see that
\begin{equation}\label{eq:wyei}
Y_i=K_i\setdif\bigcup_{j=1}^{m_i}N_{ij}
\end{equation}
since $p_1|_{L_i}$ has inverse $s\mapsto(s,\alp_i(s))$.  We also
see from (ii) that $\alp_i:K_i\to G$ is an affine map and
that $\alp_i|_{Y_i}=\alp|_{Y_i}$. \endpf

Now let us suppose that $G$ and $H$ are topological groups and the topology
of $G$ is locally compact and Hausdorff.  If $S$ is any subset of $H$
we let $\wbar{S}$ be the closure of $S$.  The following result for abelian 
groups is given in \cite[4.2.4 \& 4.5.2]{rudin}.  While our proof of 
(i) differs from the one given there, the proof of (ii) is similar
and is included for convenience of the reader.  We let $\ocringh$ denote
the {\it open coset ring}, the smallest ring of subsets of $H$ containing
all open cosets.

\begin{context}\label{lem:context}
{\bf (i)} If $\alp:C\subset H\to G$ is affine, and continuous on $C$, then it 
admits a continuous extension to an affine map $\bar{\alp}:\wbar{C}\to G$

{\bf (ii)} If $\alp:Y\subset H\to G$ is piecewise affine and continuous
on $Y$, and $Y$ is open in $H$,
then $\alp$ admits a continuous extension $\bar{\alp}:\wbar{Y}\to G$.
Moreover, $\wbar{Y}$ is open in $H$,
admits a decomposition $\wbar{Y}=\;\,\cdot\!\!\!\!\bigcup_{i=1}^nY_i$ 
as in (\ref{eq:piecwiseaff}) where each $Y_i\in\ocringh$.
\end{context}

\proof (i) Let $r\in C$.  Then the homomorphism $t\mapsto\alp(r)^{-1}\alp(rt)$
from $r^{-1}C$ to $\alp(r)^{-1}\alp(C)$ is continuous on $C^{-1}C=r^{-1}C$,
and hence left uniformly continuous by \cite[5.40(a)]{hewittrI} -- i.e.\
for every (compact) neighbourhood $W$ of $e_G$, the unit in $G$,
there is a neighbourhood $U$ of $e_H$ in $H$ such that
\begin{equation}\label{eq:unifcont}
\text{if }s^{-1}t\in U\text{ then }
\alp(rs)^{-1}\alp(rt)=\alp(rs)^{-1}\alp(r)\alp(r)^{-1}\alp(rt)\in W.
\end{equation}
If $s_0\in\wbar{C}$, then any net $(s_i)_i$ from $C$ which converges
to $s_0$ is left Cauchy: for any neighbourhood $U$ of $e_H$, there
is $i_U$ such that $i,j\geq i_U$ implies that $s_i^{-1}s_j\in U$.
Hence the net $(r^{-1}s_i)_i$ is left Cauchy as well, and so if
$U$ is chosen to satisfy (\ref{eq:unifcont}) and $i,j\geq i_U$ then
\[
\alp(s_i)^{-1}\alp(s_j)=\alp(rr^{-1}s_i)^{-1}\alp(rr^{-1}s_j)\in W
\]
so $(\alp(s_i))_i$ is left Cauchy in $G$.  However, Hausdorff locally 
compact groups are complete, whence we obtain a unique limit
$\bar{\alp}(s_0)=\lim_i\alp(s_i)$.  That $\bar{\alp}:\wbar{C}\to G$ 
is affine follows from continuity of the group operations. 

(ii) Since we assume that $\alp:Y\subset H\to C$ is piecewise affine, we
can decompose
\[
Y=\disunion{i=1}{n}Y_i\quad\text{where each}\quad 
Y_i=K_i\setdif\bigcup_{j=1}^{m_i}N_{ij}
\]
as in (\ref{eq:wyei}).  We may reorder the indices so that
$\wbar{Y}_1,\dots,\wbar{Y}_{n'}$ represents the collection
of closures having nonempty interiors.  
For $i=1,\dots,n'$ it follows that the coset $\wbar{K}_i$ must have nonempty
interior, and hence is open.  
We may reorder the second indices so that for each $i$,
$N_{i1},\dots,N_{im_i'}$ is the collection of cosets which have nonempty
interiors in $\wbar{K}_i$, and hence are both closed and open.  Then
\[
\wbar{Y}_i=\wbar{K}_i\setdif\bigcup_{j=1}^{m_i'}N_{ij}
\]
and hence is open.  Thus
\[
\wbar{Y}=\disunion{i=1}{n'}\wbar{Y}_i.
\]
Since $\alp|_{Y_i}$ is continuous, the affine map $\alp_i:K_i\to G$
such that $\alp_i|_{Y_i}=\alp|_{Y_i}$ is continuous on $Y_i$, and, by
uniformity of the topology, continuous on $K_i$.
By (i) we may extend $\alp_i:K_i\to G$ to a continuous affine map
$\bar{\alp}_i:\wbar{K}_i\to G$.
We thus let $\bar{\alp}:\wbar{Y}\to G$ be determined by
$\bar{\alp}|_{\wbar{Y}_i}=\bar{\alp}_i|_{\wbar{Y}_i}$ for each $i=1,\dots,n'$.
\endpf

We note that if $C$ is a closed coset, in the hypotheses of (i) above,
then it may not be the case that $\alp(C)$ is closed in $G$.  Consider,
for example, the map $n\mapsto e^{in}$ from the integers $\Zee$ to the circle
group $\Tee$.  See Corollary \ref{cor:mapfalg} below, for a condition which
guarantees that the range of $\alp$ is closed.

Let $G$ and $H$ be locally compact Hausdorff groups, which are usually 
referred to as simply ``locally compact''. If $\alp:Y\subset H\to G$,
then we say $\alp$ is a {\it continuous piecewise affine map} if

(i) $\alp$ is piecewise affine, and

(ii) $Y$ is both open and closed in $H$.

\section{On Idempotents in Fourier-Stieltjes Algebras}

The major result of Host \cite{host}
states that any idempotent in $\fsalg$, for any locally compact group $G$,
is the indicator function of a set from $\ocring$, the ring of sets
generated by cosets of open subgroups.  
While the following proposition borrows from Host's methods, it cannot
be directly deduced from \cite{host}.

\begin{contractiveidempotent}\label{theo:contractiveidempotent}
Let $G$ be a locally compact group and $u$ be an idempotent in $\fsalg$.

{\rm (i)} $\norm{u}=1$ if and only if $u=1_C$ for some open coset $C$ in $G$.

{\rm (ii)} $u$ is positive definite if and only if $u=1_H$ for some open 
subgroup $H$ of $G$.
\end{contractiveidempotent}

\proof Sufficiency for each of (i) and (ii) above is well known.  Let us 
review it, briefly.  If $H$ is an open subgroup of $G$ let $G/H$ denote
the discrete space of left cosets of $H$ and 
$\pi_H:G\to\fU\left(\ell^2(G/H)\right)$ the {\it quasi-left regular 
representation}, given by $\pi_H(s)\del_{tH}=\del_{stH}$.  Then
$1_H=\inprod{\pi_H(\cdot)\del_H}{\del_H}$, and hence is positive definite
with $\norm{1_H}=1_H(e)=1$.  If $C=sH$ for some $s\iin G$ then
\[
1_C=\inprod{\pi_H(\cdot)\del_H}{\del_{sH}}=
\inprod{\pi_H(s)^*\pi_H(\cdot)\del_H}{\del_H}=s\con 1_H
\]
and hence $\norm{1_C}=\norm{s\con 1_H}=\norm{1_H}=1$.  Thus it remains
to prove necessity for each (i) and (ii).

(i)  By \cite{eymard} there exist a continuous unitary representation
$\pi:G\to\fU(\fH)$ and vectors $\xi,\eta\iin\fH$ such that
\[
u=\inprod{\pi(\cdot)\xi}{\eta}\quad\aand\quad
\norm{\xi}=\norm{\eta}=\norm{u}=1.
\]
Since $\pi(s)$ is a unitary for any $s\iin G$, $\norm{\pi(s)\xi}=1$, and 
hence, by the Cauchy-Schwarz inequality we have that
\[
u(s)=\inprod{\pi(s)\xi}{\eta}=1\quad\text{ if and only if }\quad
\pi(s)\xi=\eta.
\]
Since $u$ is idempotent, we thus see that 
\[
u(s)=\inprod{\pi(s)\xi}{\eta}=
\begin{cases} 1 & \text{if }\pi(s)\xi=\eta, \\ 0 & \text{otherwise.}\end{cases}
\]
Let $C=\supp{u}=\{s\in G:u(s)=1\}$, and then we have that
\[
C=\{s\in G:\pi(s)\xi=\eta\}=\{s\in G:\xi=\pi(s)^*\eta\}.
\]
If $r,s,t\in C$, then we have that
\[
\pi(rs^{-1}t)\xi=\pi(r)\pi(s)^*\pi(t)\xi=\pi(r)\pi(s)^*\eta
=\pi(r)\xi=\eta
\]
so $rs^{-1}t\in C$.  Hence $C$ is a coset by Proposition \ref{prop:coset}.
$C$ is open since $u$ is continuous. 

(ii) First, since $\norm{u}=u(e)$, and $u$ is idempotent, we have either
that $u=0$ or $\norm{u}=1$.  The first case is trivial.  In the second, we
have from (i) that $u=1_C$ for some open coset $C$ in $G$.  Since 
$e\in C$, we must have that $C$ itself is a subgroup. \endpf

We note that if $G$ is an abelian group, then for any nontrivial idempotent
$u$ of $\fsalg$ we have either that $\norm{u}=1$ or $\norm{u}\geq
\frac{1}{2}(1+\sqrt{2})$ by \cite{saeki}.  We do not know if a similar result
holds for nonabelian groups.

The next result is well-known and a proof can be found in \cite{neumann}.
For abelian groups an interesting analytic proof is given in 
\cite[4.3.3]{rudin}.  We offer another analytic proof which is valid for
any infinite group.

\begin{fcfail}\label{prop:fcfail}
It is impossible to cover a group $G$ with finitely many cosets of subgroups 
of infinite index.
\end{fcfail}

\proof We will consider $G$ to be a discrete group.  Let $\wapg$ be the
C*-algebra of {\it weakly almost periodic functions} on $G$ 
(see \cite{burckel}).  Then $\wapg\supset\fsalg$ and $\wapg$ admits
a unique translation invariant mean $m$.  Let $\Xi$ denote the ring of 
subsets $E$ of $G$
such that $1_E\in\wapg$.  Then $\Xi\supset\cring$ by \cite{host} (or by Theorem
\ref{theo:contractiveidempotent} above).
Let $\til{m}:\Xi\to [0,1]$ be the finitely additive measure given by
$\til{m}(E)=m(1_E)$.  Then $\til{m}(G)=1$.
Thus if $C$ is a coset of a subgroup of infinite index
in $G$, $\til{m}(C)=0$.  Thus for any finite collection $C_1,\dots,C_n$
of such, $\til{m}(\bigcup_{i=1}^nC_i)\leq\sum_{i=1}^n\til{m}(C_i)=0$,
whence $C_1,\dots,C_n$ cannot cover $G$. \endpf

\section{The Main Result}

\subsection{Affine Maps Induce Completely Bounded Homomorphisms}
Let us begin with the converse of our main result.  Note that 
if $\fM$ and $\fN$ are von Neumann algebras with respective preduals 
$\fM_*$ and $\fN_*$, we say that a map $\Phi:\fM_*\to\fN_*$ 
is  {\it completely positive} if its adjoint, 
$\Psi^*:\fN\to\fM$ is completely positive.

\begin{affinecbhomo}\label{prop:affinecbhomo}
Let $G$ and $H$ be locally compact groups.  If $\alp:Y\subset H\to G$
is a continuous piecewise affine map, then $\Phi_\alp:\falg\to\fsalh$
given by
\begin{equation}\label{eq:phialp}
\Phi_\alp u(h)=\begin{cases} u(\alp(h)) & \text{if }h\in Y, \\
                                0       & \text{otherwise} \end{cases}
\end{equation}
is a completely bounded homomorphism.  
Moreover, $\Phi_\alp$ is completely contractive
if $\alp$ is affine, and completely positive if $\alp$ is a homomorphism
on an open subgroup.
\end{affinecbhomo}

\proof We will build the proof up in stages, beginning with homomorphisms, 
then moving to affine maps, and then to piecewise affine maps.

First, suppose that {\it $\alp$ is a homomorphism and $Y$ is an open 
subgroup}.  Then by \cite{eymard} or \cite[2.10]{arsac}, $\Phi_{\alp|_Y}:
\falg\to\fsaly$ is an isometric homomorphism whose adjoint
$\Phi_{\alp|_Y}^*:\wstary\to\vng$ is the $*$-homomorphism such that
$\Phi_{\alp|_Y}^*(\varpi_Y(y))=\lam_G\comp\alp(y)$ for each $y\iin Y$.
Since $Y$ is an open subgroup of $H$, $\fsaly$ injects contractively
into $\fsalh$ via the map which sends $v\iin\fsaly$ to the function
$\til{v}$, which takes the values of $v$ on $Y$ and $0$ otherwise.
This fact is well known, and follows from \cite[Prop.\ 1.2]{rieffel}, for 
example.  Note that $\{\til{v}:v\in\fsaly\}=m_Y\fsalh$ where
$m_Y:\fsalh\to\fsalh$ is multiplication by the idempotent $1_Y$.
Hence $\Phi_\alp$ is the composition of maps
\[
\begin{CD}
\falg@>{\Phi|_{\alp_Y}}>>\fsaly@>{v\mapsto\til{v}}>>m_Y\fsalh
\hookrightarrow\fsalh.
\end{CD}
\]
The adjoint of the inclusion map $m_Y\fsalh\hookrightarrow\fsalh$
is $m_Y^*$.  We have that $I=\varpi_H(e)=m_Y^*\varpi_H(e)=m_Y^*I$ and
$\cbnorm{m_Y^*}=\cbnorm{m_Y}=\norm{1_Y}=1$. Hence, by \cite[5.1.2]{effrosrB},
$m^*_Y$ is completely positive.  
The adjoint of the map $v\mapsto\til{v}$ is the $*$-isomorphism
$\theta:\wbar{\spn}^{w^*}\{\varpi_H(y):y\in Y\}\to\wstary$ such that
$\theta(\varpi_H(y))=\varpi_Y(y)$ for any $y\iin Y$.  Hence it follows that 
$\Phi_\alp^*=m_Y^*\comp\theta\comp\Phi_{\alp|_Y}^*$ is completely positive
and contractive.

Second, suppose that {\it $\alp$ is affine and $Y$ is an open coset}.  
Fix an element $h\iin Y$ and let $\beta:h^{-1}Y\to\alp(h)^{-1}\alp(Y)\subset G$
be the homomorphism from (\ref{eq:affhomo}).  Then $\Phi_\alp$ is the
composition of maps
\[
\begin{CD}
\falg@>{u\mapsto\alp(h)^{-1}\ast u}>>\falg@>{\Phi_\beta}>>\fsalh
@>{v\mapsto h\ast v}>>\fsalh
\end{CD}
\]
where for $s\con u(t)=u(s^{-1}t)$ and $h\con v(h')=v(h^{-1}h')$ for
$s,t\iin G$ and $h,h'$ in $H$.
Then for fixed $s\iin G$ [respectively, $h\iin H$], the translation operator
$u\mapsto s\con u$ on $\falg$ [$v\mapsto h\con v$ on $\fsalh$] has adjoint 
which is multiplication by the unitary $\lam_G(s)^*$ on $\vng$ [$\varpi_H(h)^*$
on $\wstarh$], which is a complete isometry.  Hence it follows that
$\Phi_\alp$ is a complete contraction.

Finally, we suppose that {\it $\alp$ is a continuous piecewise affine
map and $Y\in\ocringh$}.  Then, by Lemma \ref{lem:context} (ii), we can
write
\[
Y=\disunion{i=1}{n}Y_i\quad\text{where each}\quad
Y_i=L_i\setdif\bigcup_{j=1}^{m_i}M_{ij}
\]
and each $L_i$ and $M_{ij}$ is an open coset.  Moreover, for each $i$ 
there is a continuous affine map $\alp_i:L_i\to G$ such that $\alp_i|_{Y_i}=
\alp|_{Y_i}$.  Let $\ell^1(n)$ be the $n$-dimensional $\ell^1$-space
with contractive summing basis 
$\{\del_i:i=1,\dots,n\}$.  Let $A:\falg\to\fsalh\ptens\ell^1(n)$
be given for $u\iin\falg$ by the weighted amplification
\[
Au=\sum_{i=1}^n\Phi_{\alp_i}u\tens\del_i.
\]
Then $A$ is completely bounded with $\cbnorm{A}\leq n$.
Letting $\{\chi_i:i=1,\dots,n\}$ in $\ell^\infty(n)$ be the dual basis to
$\{\del_i:i=1,\dots,n\}$, we let the weighted diagonal multiplication map
$M:\fsalh\ptens\ell^1(n)\to\fsalh\ptens\ell^1(n)$ be given by
\[
M=\sum_{j=1}^n m_{1_{Y_i}}\tens\dpair{\cdot}{\chi_j}\del_j\quad
\text{so}\quad M\sum_{i=1}^nv_i\tens\del_i=\sum_{i=1}^n1_{Y_i}v_i\tens\del_i.
\]
Hence $M$ is completely bounded with 
$\cbnorm{M}\leq\sum_{i=1}^n\norm{1_{Y_i}}$.  Let $\trace$ be
the linear functional on $\ell^1(n)$ implemented by $\chi_1+\dots+\chi_n$,
so that $\id_{\fsalh}\tens\trace:\fsalh\ptens\ell^1(n)\to\fsalg$ is given by
\[
\id_{\fsalh}\tens\trace\left(\sum_{i=1}^nv_i\tens\del_i\right)
=v_1+\dots+v_n
\]
and $\cbnorm{\id_{\fsalh}\tens\trace}\leq\norm{\trace}=1$.  Then we
have that $\Phi_\alp$ is the composition of maps
\[
\begin{CD}
\falg@>{A}>>\fsalh\ptens\ell^1(n)@>{M}>>\fsalh\ptens\ell^1(n)
@>{\id_{\fsalh}\tens\trace}>>\fsalh
\end{CD}
\]
and is thus completely bounded. \endpf

\begin{affinecbhomo1}\label{cor:affinecbhomo1}
If $\alp:Y\subset H\to G$ is a continuous piecewise affine map
then the map $\Psi_\alp:\fsalg\to\fsalh$ given similarly as $\Phi_\alp$ in 
(\ref{eq:phialp})
is a completely bounded homomorphism.  Moreover,
$\Psi_\alp$ is completely contractive if $\alp$ is affine, and 
$\Psi_\alp$ is completely positive if $\alp$ is a homomorphism on an open 
subgroup.  If $G$ is amenable, then $\cbnorm{\Psi_\alp}=\cbnorm{\Phi_\alp}$.
\end{affinecbhomo1}

\proof That $\Psi_\alp$ has the properties claimed
can be seen by making trivial modifications to the proof of 
the proposition above: replace $\lam_G$ by $\varpi_G$, and $\vng$ by $\wstarg$.

If $G$ is amenable, there is a norm 1 bounded approximate identity
$(u_i)_i$ for $\falg$ by \cite{leptin}.  If $v\in\fsalg$ with $\norm{v}=1$,
then $(vu_i)_i$ is a contractive net which
converges to $v$ in the multiplier topology, and hence uniformly on compact 
subsets of $G$.
Hence $\Phi_\alp(vu_i)\to \Psi_\alp v$ uniformly on compact subsets
of $G$. Since $(\Phi_\alp(vu_i))_i$ is bounded by $\norm{\Phi_\alp}$
it follows that it converges to $\Psi_\alp v$ in the weak* topology, 
and hence $\norm{\Psi_\alp v}\leq\norm{\Phi_\alp}$.
Thus $\norm{\Psi_\alp}\leq\norm{\Phi_\alp}$.  

We can obtain that
$\cbnorm{\Psi_\alp}\leq\cbnorm{\Phi_\alp}$, similarly.  Indeed, for
any $n\cross n$ matrix $[v_{kl}]$ over $\fsalg$, the amplified 
multiplication maps $[v_{kl}]\mapsto [v_{kl}u_i]$ are contractive.  We get for
any $k,l$ that $\Phi_\alp(v_{kl}u_i)\to \Psi_\alp v_{kl}$ weak*, and hence
$\left[\Phi_\alp(v_{kl}u_i)\right]\to\left[\Psi_\alp v_{kl}\right]$ weak*,
from which it follows that $\norm{\Psi_\alp^{(n)}}\leq\norm{\Phi_\alp^{(n)}}$.
Hence $\cbnorm{\Psi_\alp}\leq\cbnorm{\Phi_\alp}$.
Since $\Phi_\alp=\Psi_\alp|_{\falg}$, we obtain the converse inequality,
$\cbnorm{\Phi_\alp}\leq\cbnorm{\Psi_\alp}$. \endpf

\begin{notnormext}\label{rem:notnormext}
If $G=\ftwo$, the free group on two generators, then it may not be the 
case that $\cbnorm{\Phi_\alp}=\cbnorm{\Psi_\alp}$ for a piecewise affine map
$\alp:Y\subset\ftwo\to\ftwo$ which is not itself affine.
{\rm To see this let $a$ and $b$ be the generators for $\ftwo$.  For
any $n=1,2,\dots$ the set $E_n=\{a^kb^k:k=1,\dots,n\}$ is a {\it free set}
in the sense of \cite{leinert}.  Let $\alp_n:E_n\to\ftwo$
be the inclusion map, which is affine as $E_n$ is finite.  Then
\begin{equation*}
\sup_n\cbnorm{\Phi_{\alp_n}}\leq 2\quad\text{while}\quad
\sup_n\cbnorm{\Psi_{\alp_n}}=\infty.
\end{equation*}
Indeed, first observe that for any $u\iin\fal{\ftwo}$, 
$\Phi_{\alp_n}u=1_{E_n}u$.  Then, by \cite[Bem.\ (13)]{leinert1},
$1_{E_n}\in\mathrm{B}_2(G)$, the algebra of Herz-Schur multipliers,
and is of norm no greater than 2.
Then it follows \cite{bozejkof} (or \cite{spronk1}) that
$u\mapsto 1_{E_n}u$ is completely bounded on $\falg$.
(We note that an alternative line of proof can be followed,
using \cite[Sec.\ 2]{bozejkof1} or \cite[Sec.\ 2]{pisierb},
by which we may see that $\sup_n\cbnorm{\Phi_{\alp_n}}\leq 25$.
These proofs make use of the space $\mathrm{T}_1(G)$ of Littlewood
functions.)  On the other hand, note that $\Psi_{\alp_n}1_{\ftwo}=1_{E_n}$.  
If the sequence
$\{1_{E_n}:n=1,2,\dots\}$ were bounded in $\fsalg$, then its
pointwise limit $1_E$, where $E=\{a^nb^n:n=1,2,\dots\}$,
would be in $\fsalg$.  However, this would contradict 
\cite[Lem.\ 2.7]{pisierb}.} \endpf
\end{notnormext}

\subsection{Some Lemmas}
We will need two lemmas to proceed to our main result.  

The first one
gives a construction of a nice {\it bounded approximate diagonal} for
$\falgg$.  This construction is a simplified version of the one in
the seminal paper \cite{ruan} (see also \cite[Sec.\ 7.4]{rundeB}).

\begin{bad}\label{lem:bad}
Let $G$ be an amenable locally compact group.  Then there exists in
$\falgg$ a net $(w_l)_l$ such that 

{\rm (i)} each $w_l$ is positive definite and of norm $1$, and

{\rm (ii)} for each $(s,t)\in G\cross G$
\[
\lim_l w_l(s,t)=\begin{cases} 1 &\text{if }s=t, \\ 
                              0 &\text{otherwise.} \end{cases}
\] 
\end{bad}

\proof First we will obtain a {\it positive approximate indicator}
for the diagonal subgroup as in \cite{aristovrs}.
Since $G$ is amenable, there is a positive contractive quasi-central 
approximate identity $(e_i)_i$ for $\bloneg$ by \cite{losertr}.
Moreover, by \cite{stokke}, $(e_i)_i$ can be chosen to have compact
supports tending to the identity in $G$.  We let $\xi_i=e_i^{1/2}$
for each $i$ and define for $(s,t)\iin G\cross G$
\[
u_i(s,t)=\inprod{\lam_G(s)\rho_G(t)\xi_i}{\xi_i}
\]
where $\rho_G:G\to\fU(\bltwog)$ is the right regular representation of $G$.
Then each $u_i$ is a positive definite function with $\norm{u_i}=
u(e,e)= 1$.  It follows, using computations exactly as in 
\cite[Theo.\ 2.4]{aristovrs}, that 
\[
\lim_i u_i(s,t)=\begin{cases} 1 &\text{if }s=t, \\ 
                              0 &\text{otherwise.} \end{cases}
\] 
(We note that if $G$ is a small invariant neighbourhood group, there
is a very simple construction of a net, having the properties of 
$(u_i)_i$, in \cite{spronk}.)

Since $G$ is amenable, $G\cross G$ is amenable and hence
there exists a net $(v_j)_j$ of $\falgg$
comprised of norm 1 positive definite functions which converges  
uniformly on compact sets to $1_{G\times G}$.
To see this, by \cite[Prop.\ 6.1]{hulanicki} (also see \cite[4.21]{paterson})
there exists a bounded
net of elements $(v'_j)_j$ from $\falgg^+$ which converges  
uniformly on compact sets to $1_{G\times G}$, and normalise
by taking $v_j=\frac{1}{v'_j(e,e)}v'_j$.  (Note that
by \cite{granirerl}, this net is a bounded approximate identity
for $\falgg$. This gives an alternative proof to \cite{leptin}.)  
Then the product net with elements
\[
w_l=v_j u_i
\]
is the desired bounded approximate diagonal.  \endpf

The next lemma is a straightforward technical result on complete positivity.
Though it is surely well-known, we include a proof for convenience of the
reader.

\begin{cpmap}\label{lem:cpmap}
Let $\fM_1,\fM_2$ and $\fN$ be von Neumann algebras and let 
\linebreak $\Phi:(\fM_1)_*\to(\fM_2)_*$  be a
completely positive map.  Then $\Phi\tens\id_{\fN_*}:
(\fM_1)_*\ptens \fN_*$ $\to(\fM_2)_*\ptens\fN_*$ is also
completely positive.
\end{cpmap}

\proof By definition we have that $\Phi^*:\fM_2\to\fM_1$ 
is completely positive.  Hence the normal
map $\Phi^*\tens\id_\fN:\fM_2\vntens\fN\to\fM_1\vntens\fN$ makes sense.
Moreover, it is simple to verify, using elementary tensors in 
$\fM_2\vntens\fN$, that $\Phi^*\tens\id_N=(\Phi\tens\id_{\fN_*})^*$,
where we identify 
$\fM_i\vntens\fN\cong\left((\fM_i)_*\ptens\fN_*\right)^*$ for $i=1,2$.
Thus we need only to see that $\Phi^*\tens\id_\fN$ is completely positive.

First, let $b=\Phi^*(I)^{1/2}$.  By \cite[5.1.6]{effrosrB} there is a unital
completely positive map $\theta:\fM_2\to\fM_1$ such that
$\Phi^*=b\theta(\cdot)b$.  Hence
\[
\Phi^*\tens\id_\fN=(b\tens I)\theta\tens\id_\fN(\cdot)(b\tens I).
\]
Now $\theta\tens\id_\fN$ is completely positive since 
$\theta\tens\id_\fN(I\tens I)=I\tens I$ while $\cbnorm{\theta\tens\id_\fN}
\leq\cbnorm{\theta}=\norm{\theta(I)}=1$, and we employ \cite[5.1.2]{effrosrB}.
Thus $\Phi^*\tens\id_\fN$ is a composition of completely positive maps
whence it is completely positive.  \endpf

\begin{inversemap}\label{rem:inversemap}{\rm
If $G$ is any locally compact group, the map $\iota:\falg\to\falg$
given by $\iota u(s)=u(s^{-1})$ ($s\iin G$) is positive -- it takes
positive definite elements to positive definite elements --
contractive and a homomorphism.  
However, by \cite[Prop.\ 1.5]{forrestr}, 
$\iota$ is
completely bounded only if $G$ has an abelian subgroup of finite index.
Thus if $G$ does not admit an abelian subgroup of finite index, then
$\falg$ admits a bounded homomorphism which is not completely bounded.}
\end{inversemap}

\subsection{The Main Result}

\begin{maintheorem}\label{theo:maintheorem}
Let $G$ and $H$ be locally compact groups with $G$ amenable, and let
$\Phi:\falg\to\fsalh$ be a completely bounded homomorphism.  Then there
is a continuous piecewise affine map $\alp:Y\subset H\to G$ such that for each
$h\iin H$
\[
\Phi u(h)=\begin{cases} u(\alp(h)) & \text{if } h\in Y, \\
                        0   & \text{otherwise.}  \end{cases}
\]
Moreover, $\alp$ is affine if $\Phi$ is completely contractive, and 
$\alp$ is a homomorphism defined on an open subgroup if $\Phi$ is 
completely positive.
\end{maintheorem}

\proof First, it will be convenient for us to let 
$G_\infty=G\dot{\cup}\{\infty\}$ be either the one point compactification of 
$G$ if $G$ is not compact, or the topological coproduct
if $G$ is compact.  Each element $u$ of $\falg$ extends
to a continuous function on $G_\infty$ by setting $u(\infty)=0$.
Now, as observed in \cite{ilie}, since $G$ is the Gel'fand spectrum
of $\falg$, for each $h\iin H$ there is an $\alp(h)\iin G_\infty$
such that
\begin{equation}\label{eq:alpdef}
\Phi u(h)=u(\alp(h)).
\end{equation}
%
The map $\alp:H\to G_\infty$ is continuous.  Indeed, suppose not.
Then there is an $h_0$ in $H$, a neighbourhood $U$ of $\alp(h_0)$ and a net
$(h_i)_i$ in $H$ such that $h_i\to h_0$ but $\alp(h_i)\not\in U$ for any $i$.
If $\alp(h_0)\in G$, find a $u$ in $\falg$ such that $\supp{u}\subset U$
and $u(\alp(h_0))=1$.  Then
\[
1=u(\alp(h_0))=\Phi u(h_0)=\lim_i \Phi u(h_i)=\lim_i u(\alp(h_i))=0
\]
which is absurd.  If $\alp(h_0)=\infty$, then $K=G_\infty\setdif U$ is 
compact.
Find $u\iin \falg$ such that $u|_K=1$ and $\supp{u}$ is compact, and we
reach a similar contradiction as above.  We thus have that
\begin{equation}\label{eq:alpcont}
\alp:H\to G_\infty\text{ is continuous and }Y=\alp^{-1}(G)\text{ is open in }H.
\end{equation}
Note that
\begin{equation}\label{eq:wyechar}
Y=\{h\in H:\text{ there exists a }u\iin \falg\text{ such that }
\Phi u(h)\not=0\}.
\end{equation}

We will now suppose that $\Phi$ is completely bounded, and address the 
cases that it is completely contractive or completely positive later.
Then the map \linebreak
$\Phi\tens\id_{\falg}:\falg\ptens\falg\to\fsalh\ptens\falg$
is completely bounded with \linebreak
$\cbnorm{\Phi\tens\id_{\falg}}\leq\cbnorm{\Phi}$.  We can identify
$\falg\ptens\falg$ completely isometrically with $\falgg$ via
\begin{equation}\label{eq:tenstocross}
u\tens v\mapsto u\cross v
\end{equation}
where $u\cross v(s,t)=u(s)v(t)$ for $(s,t)\iin G\cross G$.
Indeed, the adjoint of this map is the isomorphism $\vngg\cong\vng\vntens\vng$,
which is spatially implemented.
We can also identify $\fsalh\ptens\falg$ completely isometrically
as a subspace of $\fsalhcg$ via a map like (\ref{eq:tenstocross}), where
$u\cross v(h,t)=u(h)v(t)$ for $(h,t)\iin H\cross G$.  Indeed, the
adjoint of such a map would be the canonical homomorphism from
$\wstarhcg$ to $\wstarh\vntens\vng$, which extends 
$\vpi_{H\times G}(h,t)$ $\mapsto\vpi_H(h)\tens\lam_G(t)$.  
This is a complete quotient map.  
Furthermore, the inclusion $\fsalhcg\hookrightarrow\fsalhcgd$ is
a completely positive isometry.  Thus the composition of maps
\[
\begin{CD}
\falgg\cong\falg\ptens\falg@>{\Phi\tens\id_{\falg}}>>
\fsalh\ptens\falg\hookrightarrow\fsalhcg\hookrightarrow\fsalhcgd
\end{CD}
\]
forms a map $J_\Phi$, such that $\norm{J_\Phi}\leq\cbnorm{J_\Phi}\leq
\cbnorm{\Phi}$.  On any elementary product of functions $u\cross v$ 
in $\falgg$ and for any $(h,t)\iin H\cross G$ we have
\[
J_\Phi u\cross v(h,t)=(\Phi u)\cross v (h,t)=u(\alp(h))v(t)
\]
where we recall that $u(\infty)=0$.  We can extend this to linear
combinations of elementary products and hence, by continuity, to any element 
$w$ in $\falgg$, so we obtain for any $(h,t)\iin H\cross G$
\begin{equation}\label{eq:jphiform}
J_\Phi w(h,t)=\begin{cases} w(\alp(h),t) & \text{if }h\in Y, \\
                               0 & \text{otherwise.} \end{cases}
\end{equation}
Now, let $(w_l)_l$ be the net from Lemma \ref{lem:bad}, and consider the net
$\left(J_\Phi w_l\right)_l$ in $\fsalhcgd$.  We have for each 
$(h,t)\iin H\cross G$ that
\begin{align*}
\lim_l J_\Phi w_l(h,t)
&={\begin{cases} \lim_l w_l(\alp(h),t) & \text{if }h\in Y, \\
                               0 & \text{otherwise} \end{cases}} \\
&={\begin{cases} 1 & \text{if }h\in Y\aand\alp(h)=t,\\
                               0 & \text{otherwise.} \end{cases}}
\end{align*}
Thus $\left(J_\Phi w_l\right)_l$ converges pointwise to $1_{\Gamma_\alp}$, 
where $\Gamma_\alp$ is the graph of $\alp$, as in (\ref{eq:graph}).  
Since $\left(J_\Phi w_l\right)_l$ is bounded, we
conclude that $1_{\Gamma_\alp}\in\fsalhcgd$, which implies by \cite{host} that
$\Gamma_\alp\in\cringhcg$.  
This implies, by Lemma \ref{lem:graphs} (iii), that $\alp$ is piecewise
affine.  Since (\ref{eq:alpcont}) holds, we can apply Lemma \ref{lem:context}
to see that $\alp$ extends to a continuous piecewise affine map
$\bar{\alp}:\wbar{Y}\to G$.  However, it follows from this that $Y$ is closed.
Indeed if $h_0\in\wbar{Y}$, there is an element $u$ of $\falg$ such that
$u(\bar{\alp}(h_0))\not=0$.  Then for any net $(h_i)_i$ in $Y$, converging
to $h_0$ we have
\[
\Phi u(h_0)=\lim_i \Phi u(h_i)=\lim_i u(\alp(h_i))=u(\bar{\alp}(h_0))\not=0
\]
and it follows from (\ref{eq:wyechar}) that $h_0\in Y$.
Hence $\bar{\alp}=\alp$ and $\alp$ is itself a continuous piecewise affine map.

Now suppose that $\Phi$ is completely contractive.  It then follows that
the net $\left(J_\Phi w_l\right)_l$ is completely contractive, and hence
$1_{\Gamma_\alp}\in\fsalhcgd$ with $\norm{1_{\Gamma_\alp}}\leq 1$.  By
Theorem \ref{theo:contractiveidempotent} we conclude that $\Gamma_\alp$
is a coset, which in turn forces $\alp$ to be an affine map by
Lemma \ref{lem:graphs} (ii).

Finally, if $\Phi$ is completely positive, then $\Phi\tens\id_{\falg}$ is a 
completely positive map by Lemma \ref{lem:cpmap}, and hence $J_\Phi$ is 
completely positive.
Thus the net $\left(J_\Phi w_l\right)_l$ is a net of norm
1 positive definite functions.  It follows then that $1_{\Gamma_\alp}$ is
positive definite, which, by Theorem \ref{theo:contractiveidempotent} 
forces $\Gamma_\alp$ to be a group.  Hence
$\alp$ is a homomorphism on some open subgroup by Lemma \ref{lem:graphs} (i). 
\endpf

The theorem stated above appears to be the best result possible.
We have already seen in Remark \ref{rem:inversemap} that there
exist bounded homomorphisms of Fourier algebras which are not completely
bounded.  Moreover,
we can show for a large class of nonamenable groups, namely those which
contain a closed noncommutative free group, that Theorem
\ref{theo:maintheorem} fails.  We note that a {\it nonamenable} 
group $G$ contains
a discrete copy of a noncommutative free group if $G$ is almost connected
(see \cite[3.8]{paterson}), or if $G$ is linear with the discrete topology
(see \cite{tits} or, for some special cases, see \cite{harpe}).  We are
indebted to B.E.\ Forrest for indicating the following example to us.

\begin{failure}\label{prop:failure}
If $G$ is a nonamenable group which contains a discrete copy of the free
group $\ftwo$, then there exists a completely bounded homomorphism
$\Phi:\falg\to\falftwo$ which is not of the form $\Phi_\alp$ as in
(\ref{eq:phialp}) for a piecewise affine map $\alp$.
\end{failure}

\proof The restriction map $u\mapsto u|_{\ftwo}$ from $\falg$ to
$\falftwo$ is a contractive quotient map by \cite{herz}.  Moreover, the
adjoint of this restriction map is the $*$-homomorphism from
$\vnftwo$ to $\vng$ which extends $\lam_{\ftwo}(s)\mapsto\lam_G(s)$,
so it is completely positive.  Now let $a$ and $b$ be the generators
for $\ftwo$ and $E=\{a^nb^n:n=1,2,\dots\}$, so $E$ is a {\it free set}.
As in Remark \ref{rem:notnormext},
we see that the map $v\mapsto 1_Ev$ is a completely bounded map
on $\falftwo$.  We let $\Phi$ be the composition of maps
\[
\begin{CD}
\falg@>{u\mapsto u|_{\ftwo}}>>\falftwo@>{v\mapsto 1_Ev}>>\falftwo.
\end{CD}
\]
Then $\Phi=\Phi_\alp$ where $\alp:E\subset\ftwo\to G$
is the inclusion map.  However  $1_E\not\in\fsal{\ftwo}$ by
\cite[Lem.\ 2.7]{pisierb}.  Hence, by \cite{host},
$E\not\in\cringftwo$, so $\alp$ is not piecewise affine.  \endpf

\subsection{Some Consequences}
In what follows we will always let $G$ and $H$ be locally compact groups 
with $G$ amenable.

The first one is a direct application of Theorem \ref{theo:maintheorem}
and Corollary \ref{cor:affinecbhomo1}.

\begin{maintheorem1}
Any completely bounded homomorphism $\Phi:\falg\to\fsalh$ extends to a
completely bounded homomorphism $\Psi:\fsalg\to\fsalh$
with $\cbnorm{\Psi}=\cbnorm{\Phi}$.
\end{maintheorem1}

The next corollary follows from the fact that a connected group
admits no continuous piecewise affine maps which are not themselves affine.
See Lemma \ref{lem:context}.

\begin{hconnected}
If $H$ is connected, then any completely bounded homomorphism
$\Phi:\falg\to\fsalh$ is of the form (\ref{eq:phialp}) for an affine map
$\alp$.  In particular, $\Phi$ is completely contractive.
\end{hconnected}

Now we will consider homomorphisms between $\falg$ and $\falh$. 
If $Y$ and $X$ are locally compact spaces, we say that a map $\alp:Y\to X$
is {\it proper} if $\alp^{-1}(K)$ is compact in $Y$ for every compact 
subset $K$ of $X$. For abelian groups the result below can be found in 
\cite{cohen}.

\begin{mapfalg}\label{cor:mapfalg}
A map $\Phi:\falg\to\falh$ is a completely bounded homomorphism if and only if
$\Phi=\Phi_\alp$ as in (\ref{eq:phialp}) where $\alp$ is a continuous
piecewise affine map which is proper.
In particular, $\alp$ is a closed map.
\end{mapfalg}

\proof As in the proof of Theorem \ref{theo:maintheorem}, let us employ
the convention that $\alp:H\to G_\infty$.

Since $\Phi:\falg\to\falh\subset\fsalh$, the existence of
a piecewise affine map $\alp$ such that $\Phi=\Phi_\alp$ follows from
our main result.  Now if there were a compact subset $K$ of $G$
such that $\alp^{-1}(K)$ is not compact, then for any $u\in\falg$
such that $u|_K=1$, we would have that $\Phi u=u\comp\alp$ would not
vanish at infinity, and hence would not be in $\falh$.  Hence necessity
is proven.

To obtain the converse, 
let $\falcg$ denote the subspace of $\fsalg$ consisting of
functions of compact support.  Then $\falcg$ is a dense subspace of $\falg$.
Similarly define $\falch$.  If
$\Phi=\Phi_\alp$ for a continuous piecewise affine map $\alp$, then
$\Phi:\falg\to\fsalh$ is a completely bounded homomorphism by Proposition
\ref{prop:affinecbhomo}.  Since for any $u\iin\falcg$, $\supp{u\comp\alp}
=\alp^{-1}(\supp{u})$, we have that $\Phi(\falcg)\subset\falch$, and hence
$\Phi(\falg)\subset\falh$. 

Now we shall see that $\alp$ is closed.
Let $E$ be a non-empty closed subset of $H$,
$s_0\in\wbar{\alp(E)}$ and $\fU$ be a neighbourhood basis
at $s_0$ consisting of relatively compact sets.  Then
$\{\alp^{-1}(\wbar{U})\cap E\}_{U\in\fU}$ is a family of compact subsets of $H$
having the finite intersection property, whence 
$L=\bigcap_{U\in\fU}\alp^{-1}(\wbar{U})\cap E\not=\varnothing$.
If $h_0\in L$, then $\alp(h_0)=s_0$, since $\alp(h_0)$ is contained in
$\wbar{U}$ for each $U\iin\fU$.
\endpf

We can thus obtain an analogue of Walter's Theorem \cite{walter}.

\begin{walter}
Let $\Phi:\falg\to\falh$ be an completely contractive isomorphism.
Then there exists an element $s_0\iin G$ and a topological group isomorphism 
$\beta:H\to G$ such that
\[
\Phi u(h)= u(s_0\beta(h))
\]
for each $h\iin H$.  Moreover, if $\Phi$ is completely positive then $s_0=e$.
\end{walter}

\proof Let $\alp:Y\subset H\to G$ be the affine map whose existence
is guaranteed by Corollary \ref{cor:mapfalg}.  Since $\Phi$ is surjective
and $\falh$ is a point-separating regular algebra on $H$,
$Y=H$ and $\alp$ is injective.  Since $\Phi$ is injective, $\alp(H)$ is 
dense in $G$.  Hence, as $\alp$ is closed, we obtain that $\alp(H)=G$.  
Thus $\alp$ is bijective, and hence open, so it is a 
homeomorphism.  Let $s_0=\alp(e_H)$  and $\beta=s_0^{-1}\alp(\cdot)$,
so $\beta$ is a homomorphism as in (\ref{eq:affhomo}).   

If $\Phi$ is completely positive then $\alp$ is a homomorphism
and hence $\beta=\alp$.  \endpf

It is known, to us through Z.-J.\ Ruan, that the result
above obtains in the case that
$\Phi$ is a complete isometry, without assuming that $G$ is amenable.
The proof follows the one given by Walter.

Let us close with a characterisation of the range of a completely bounded
homomorphism of Fourier algebras, which is due to Kepert \cite{kepert}
in the case that $G$ and $H$ are abelian, and due to the first author
\cite{ilie} in the case that either $G$ is discrete and amenable
or $G$ is abelian, with 
general $H$.  Below, as with all results in this section, we assume the $G$ is
amenable and $H$ is a general locally compact group.

\begin{image}
If $\Phi:\falg\to\falh$ is a completely bounded homomorphism 
with adjoint map $\Phi^*:\vnh\to\vng$ then
\[
\Phi(\falg)=\left\{u\in\falh:
\begin{matrix}  u(h)=0 \text{ if }\Phi^*(\lam_H(h))=0,\text{ and} \\
            u(h)=u(h') \text{ if }\Phi^*(\lam_H(h))=\Phi^*(\lam_H(h'))
\end{matrix}\right\}
\]
\end{image}

That $\Phi(\falg)$ is contained in the set of the latter description
is clear.  Note that, in fact, 
$\Phi^*(\lam_H(h))=\lam_G(\alp(h))$ if $ \Phi^*(\lam_H(h))\not=0$,
where $\alp$ is the proper completely affine map promised by Corollary
\ref{cor:mapfalg}.  The converse inclusion is not trivial, and its proof
is in \cite[Sec.\ 5]{ilie}.  To adapt that proof to be sufficiently general we 
need only note Herz's extension theorem \cite{herz} --
if $F$ is a closed subgroup of
$G$ then for each $u$ in $\fal{F}$ there is an element $\til{u}\iin\falg$ such
that $\til{u}|_F=u$ and $\norm{\til{u}}=\norm{u}$; and we must use
Corollary \ref{cor:mapfalg} in place of the main theorem of \cite{ilie}.

{\bf Acknowledgments.}
We are grateful to R.R.\ Smith and B.E.\ Forrest for comments on our
presentation, and to G.\ Pisier and M.\ Bo\.{z}ejko
for clarifying comments on free sets.
We are also grateful to R.\ Stokke for pointing out an error in an earlier
version.

{
\bibliography{cbhomofalgbib}
\bibliographystyle{plain}
}

\end{document}